\documentclass[11pt,reqno]{amsart}

\usepackage{amsmath,amssymb,bbm,graphicx,enumitem,mathabx,psfrag,upgreek,xcolor}

\DeclareMathAlphabet{\matheu}{U}{eur}{m}{n}

\makeatletter
\def\@seccntformat#1{%
  \protect\textup{\protect\@secnumfont
    \ifnum\pdfstrcmp{subsection}{#1}=0 \bfseries\fi
    \csname the#1\endcsname
    \protect\@secnumpunct
  }%
}
\makeatother

\numberwithin{equation}{section}

\renewcommand{\thefigure}{\arabic{section}.\arabic{figure}}
\renewcommand{\thefigure}{\ifnum\value{section}>0 \arabic{section}.\fi\arabic{figure}}
\renewcommand{\theequation}{\ifnum\value{section}>-1 \arabic{section}.\fi\arabic{equation}}

\newtheorem{thm}[equation]{Theorem}

\newtheorem{prp}[equation]{Proposition}

\theoremstyle{definition}
\newtheorem{dfn}[equation]{Definition}
\newtheorem{ex}[equation]{Example}

\theoremstyle{remark}
\newtheorem{rem}[equation]{Remark}

\newcommand{\thmref}[1]{Theorem~\ref{#1}}

\newcommand{\dfnref}[1]{Definition~\ref{#1}}
\newcommand{\remref}[1]{Remark~\ref{#1}}
\newcommand{\exref}[1]{Example~\ref{#1}}
\newcommand{\figref}[1]{Figure~\ref{#1}}
\newcommand{\secref}[1]{Section~\textbf{\ref{#1}}}

\makeatletter
\@namedef{subjclassname@2020}{%
  \textup{2020} Mathematics Subject Classification}
\makeatother

\newlength{\mqheight}
\newlength{\mqnormalheight}
\newcommand{\mathquote}[1]{%
    \settoheight{\mqheight}{\hbox{\ensuremath{#1}}}%
    \addtolength{\mqheight}{-\mqnormalheight}%
    \text{\raisebox{\the\mqheight}{``}}%
    #1%
    \text{\raisebox{\the\mqheight}{''}}%
}

\newcommand\1{\mathbbm{1}}

\newcommand\Ac{\mathcal A}
\newcommand{\acts}{\lefttorightarrow}
\newcommand\al{\alpha}

\newcommand\De{\Delta}
\newcommand\de{\delta}

\newcommand\Fc{\mathcal F}

\newcommand\ga{\gamma}
\newcommand\gb{{\boldsymbol{g}}}

\newcommand\gr{\operatorname{gr}}
\newcommand{\Gu}{\mathbf{G}}
\newcommand{\gu}{\mathbf{g}}

\newcommand\la{\lambda}
\newcommand\lan{\langle}
\newcommand\lt{\leadsto}

\newcommand{\mapstoto}{\mathop{\,\sim\joinrel\rightsquigarrow\,}}
\newcommand\me{\text{\large$\boldsymbol{\matheu{m}}$}}

\newcommand\NN{\mathbb N}

\newcommand\nub{\text{\large$\boldsymbol{\nu}$}}

\newcommand\ol{\overline}
\newcommand\om{\omega}

\newcommand\p{\partial}
\newcommand\Pb{\mathbf P}

\newcommand\ph{\varphi}

\newcommand\ran{\rangle}
\newcommand\rb{\mathbf{r}}
\newcommand\RR{\mathbb R}
\newcommand\RW{\operatorname{RW}}

\newcommand\Sc{\mathcal S}

\newcommand{\si}{\sigma}
\newcommand{\sm}{\,\setminus\,}
\newcommand\supp{\operatorname{supp}}

\renewcommand\th{\theta}

\newcommand{\Upnu}{\matheu{N}}
\newcommand{\Upsi}{\Upsigma}

\newcommand\wh{\widehat}
\newcommand\wt{\widetilde}

\newcommand\ZZ{\mathbb Z}

\begin{document}

\title[Quasi-product Markov measures]{Quasi-product and multiplicative Markov\\ measures, Radon -- Nikodym problem,\\ random walks: the free group case}

\author{Vadim A. Kaimanovich}

\address{Department of Mathematics and Statistics, University of
Ottawa, Ottawa ON, K1N 6N5, Canada}

\email{{\tt vkaimano@uottawa.ca, vadim.kaimanovich@gmail.com}}



\begin{abstract}
The purpose of this paper is a conceptual clarification of the role of a special class of measures on topological Markov chains\,---\,multiplicative Markov measures\,---\,in the context of random walks. It is achieved by combining the interpretation of multiplicative Markov measures on the boundary of a free group as the unique solutions of the Radon -- Nikodym problem for product cocycles of the boundary action and the fact that the Radon -- Nikodym cocycle of the harmonic measure of a nearest neighbour random walk on a free group is a product one.
\end{abstract}

\maketitle

\thispagestyle{empty}

\dedicatory{\hfill \textit {To my teacher Anatoly Moiseevich Vershik,\\ \hspace*{\fill} on the occasion of his 90th birthday}}

\section*{Introduction}

{\bf 1.} The central notion of this paper is that of a \emph{quasi-product Markov measure}. It naturally arises in the presence of a \emph{topological Markov chain} structure on a finite state space (this structure amounts to specifying which transitions between states are \emph{admissible}). Then just a single \emph{transition measure} on the state space is sufficient to define a family of transition probabilities as the normalized restrictions of the transition measure to the corresponding sets of states (the transitions to which from a given departure state are admissible). The arising Markov measures on the path space (i.e., on the topological Markov chain) are called quasi-product. \emph{Multiplicative Markov measures} are those quasi-product ones for which the transition measure (called in this case the \emph{base measure}) serves not only to define the transition probabilities, but also as the initial distribution.

Quasi-product Markov measures appear already in \textsc{Parry}'s seminal 1964 paper \cite[Theorem~2]{Parry64} under the name of \emph{piecewise linear processes} (it is in this very paper that topological Markov chains were first introduced, see \textsc{Katok -- Hasselblatt} \cite[p.~744]{Katok-Hasselblatt95} for a historical and terminological discussion concerning the latter notion). They were reintroduced several times by a number of authors apparently unaware of \cite{Parry64} and of each other:

\medskip

\noindent
$\bullet\,$ \textsc{Evans} \cite{Evans85} defined \emph{quasi-product measures} (it is his terminology that we adopt) in the context of time inhomogeneous topological Markov chains (nowadays known in the theory of operator algebras as \emph{Bratteli diagrams}, e.g., see \textsc{Herman~-- Putnam~-- Skau} \cite{Herman-Putnam-Skau92}); in the time homogeneous case ($\equiv$ for \emph{stationary} Bratteli diagrams) they are precisely the quasi-product Markov measures in our sense.

\medskip

\noindent
$\bullet\,$ Later quasi-product Markov measures appeared (in a somewhat implicit form though) as the Gibbs measures of topological Markov chains whose \emph{potential depends on a single coordinate} in \textsc{Petersen~-- Schmidt} \cite[Theorem 2]{Petersen-Schmidt97} (and these are precisely the \emph{measures that are invariant under both the shift and the adic transformation}).

\medskip

\noindent
$\bullet\,$ Finally, in what concerns \emph{random walks on groups} and their \emph{harmonic measures,} the fact that the harmonic measure of a nearest neighbour random walk on a free group is Markov is essentially contained already in the pioneering work of Dynkin and his collaborators in the 60s (\textsc{Dynkin~-- \linebreak Malyutov} \cite{Dynkin-Malutov61}, \textsc{Dynkin~-- Yushkevich} \cite[pp.~142-145]{Dynkin-Yushkevich69}), and was made explicit by \textsc{Levit~-- Mol\-cha\-nov}~\cite{Levit-Molchanov71}. Its multiplicativity in the sense of the above definition was perspicuated by \textsc{Woess} \cite[Definition on p.~143]{Woess85} and \textsc{Gerl -- Woess} \cite[Definition on p.~351]{Gerl-Woess86} (who explicitly used the term \emph{multiplicative measure}). Without singling out this property terminologically, it also appears in \textsc{Ledrappier} \cite[formula~(2.7)]{Ledrappier01}. Later \textsc{Mairesse}~\cite{Mairesse05} and \textsc{Mairesse -- Matheus} \cite{Mairesse-Matheus04,Mairesse-Matheus07} reintroduced multiplicative Markov measures under the name of \emph{Markovian multiplicative} and showed that the harmonic measure on the space of infinite irreducible words belongs to this class for the nearest neighbour random walks on what they called \emph{$0$-automatic Cayley graphs} (these are the free products of cyclic or finite groups endowed with generating sets which may be larger than the ``default'' generating sets in free products). However, multiplicativity appears in all these papers \emph{post hoc}, as a consequence of an explicit description of harmonic measures.

\bigskip

{\bf 2.} In the present article we focus on the topological Markov chain determined by a \emph{finitely generated free group} $\Fc$ by fixing its free generating set and then considering the \emph{alphabet}~$\Ac$ that consists of these generators and their inverses. The admissibility condition then means that mutually inverse generators are not allowed to stand next to each other, so that the group $\Fc$ itself is identified with the set of finite admissible (\emph{irreducible,} in the group parlance) words in the alphabet $\Ac$, whereas the set of (one-sided) admissible infinite words is the group \emph{boundary} $\p\Fc$.

First we interpret multiplicative Markov measures on the boundary $\p\Fc$ as the unique solutions of the Radon -- Nikodym problem for what we call \emph{product cocycles} of the \emph{boundary action} $\Fc\acts\p\Fc$ of the free group. These are, geometrically, the simplest multiplicative real-valued positive cocycles of the boundary action. A cocycle from this class is determined by a family of positive numbers $\{\pi_a\}_{a\in\Ac}$ indexed by the alphabet~$\Ac$, which provides an assignment of weights to the (oriented) edges of the Cayley tree of the group~$\Fc$. The cocycle is then obtained by ``regularizing'' the ratios of infinite products of weights along geodesic rays of the Cayley tree by using the fact that if two such rays converge to the same boundary point, then they eventually merge.

Given a cocycle of a group action, the \emph{Radon -- Nikodym problem} for this cocycle consists in finding quasi-invariant measures on the action space whose Radon -- Nikodym cocycle coincides with the given one. It turns out that for product cocycles the Radon~-- Nikodym problem is solvable if and only if the weights $\{\pi_a\}_{a\in\Ac}$ satisfy a certain algebraic \emph{consistency condition}, in which case the solution is a unique multiplicative Markov measure, and, conversely, the Radon -- Nikodym cocycle of any such measure is a product cocycle (\thmref{thm:final} and \thmref{thm:qi}). The arising one-to-one correspondence between the base distributions $\{\nu_a\}_{a\in\Ac}$ of multiplicative Markov measures on the boundary~$\p\Fc$ (subject to the normalization condition $\sum_{a\in\Ac} \nu_a=1$) and the weights $\{\pi_a\}_{a\in\Ac}$ of the associated Radon~-- Nikodym cocycles (subject to the aforementioned consistency condition) can be easily described explicitly. The consistency condition on the weights $\{\pi_a\}_{a\in\Ac}$ also has a ``spectral'' interpretation: it is equivalent to requiring that the Perron~-- Frobenius eigenvalue of the weighted admissibility matrix be $1$ (\thmref{thm:pf}).

Thus, the class of multiplicative Markov measures on the boundary $\p\Fc$ coincides with the class of quasi-invariant probability measures with a product Radon -- Nikodym cocycle.

\bigskip

{\bf 3.} Let us add another character to the plot. This is a probability measure~$\mu$ on the alphabet $\Ac$, i.e., a collection of positive weights $\mu_a=\mu(a)$ subject to the normalization condition $\sum_{a\in\Ac}\mu_a=1$. The transitions of the associated \emph{random walk} consist in the right multiplication by $\mu$-distributed increments (geometrically, it amounts to passing to a random neighbour in the Cayley tree of the group), and its sample paths issued from the group identity are the products $g_t=h_1h_2\dots h_t$ of independent increments. As we have mentioned earlier, it was established already by \textsc{Dynkin~--} \linebreak \textsc{Malyutov}~\cite{Dynkin-Malutov61} that the sample paths of the random walk $(\Fc,\mu)$ almost surely converge to the boundary $\p\Fc$, and the resulting limit distribution $\nub$ is called the \emph{harmonic ($\equiv$ hitting) measure} of the random walk.

Now, since the random walk $(\Fc,\mu)$ is nearest neighbour, it is virtually obvious that the Radon~-- Nikodym cocycle of the harmonic measure $\nub$ is a product one, and the weights $\{\pi_a\}_{a\in\Ac}$ that determine this cocycle are the \emph{passage probabilities} of the random walk (i.e., $\pi_a$ is the probability that the random walk issued from the group identity ever visits the point $a$). Indeed, if $\Sc_a\subset\p\Fc$ denotes the ``shadow'' of a point $a\in\Ac$ cast by a light source positioned at the group identity ($\equiv$ the set of infinite irreducible words beginning with the letter $a$), then for any measurable subset ${X\subset\Sc_a}$ its harmonic measure $\nub(X)$ is equal to the product of $\pi_a$ (as one can not reach $\Sc_a$ from the group identity without passing through $a$) and the probability $a\nub(X)=\nub(a^{-1}X)$ of hitting $X$ starting from $a$, whence the Radon~-- Nikodym derivative $da\nub/d\nub$ is identically $1/\pi_a$ on $\Sc_a$.

From this point of view the earlier results on the description of the harmonic measure for a nearest neighbour random walk on the free group $\Fc$ (see above) can be reformulated in the following form: the Radon -- Nikodym cocycle of the harmonic measure $\nub$ of any nearest neighbour random walk on~$\Fc$ is a product one, and therefore the harmonic measure itself is multiplicative Markov (\thmref{thm:h}). Moreover (\thmref{thm:bij}), there is a one-to-one correspondence between three collections of numbers from the unit interval indexed by the alphabet $\Ac$:
\begin{enumerate}[label = {$\bullet$}]
\item
the weights $\mu_a=\mu(a)$ of a probability measure $\mu$ on $\Ac$ that serves as the step distribution of the random walk,
\item
the weights $\nu_a=\nub(\Sc_a)$ of the base distribution of the harmonic measure $\nub$,
\item
the passage probabilities $\pi_a$ of the random walk that also serve as the weights of the product Radon~-- Nikodym cocycle of the harmonic measure.
\end{enumerate}
The collections $\{\mu_a\}$ and $\{\nu_a\}$ are subject to the normalization conditions $\sum_{a\in\Ac}\mu_a=\sum_{a\in\Ac}\nu_a=1$, whereas the last one $\{\pi_a\}$ is subject to the aforementioned consistency condition.

An obvious consequence of \thmref{thm:bij} is that for a product cocycle~$c$ solvability of the Radon -- Nikodym problem is equivalent to its \emph{harmonicity}, i.e., to existence of a probability measure $\mu$ on the group such that $\sum_g \mu(g) c(g,\ga)=1$ for all boundary points $\ga\in\p\Fc$. It might be interesting to explore this relationship for other cocycles of the boundary action, see \remref{rem:har}.

\bigskip

{\bf 4.} All our considerations are also applicable to the Mairesse -- Matheus \mbox{0-automatic} setup, and after an appropriate reformulation most of them make sense for general topological Markov chains as well, or even for Bratteli diagrams (in particular, for general trees, see \remref{rem:cone}). However, for the sake of simplicity we focus just on the free group example, leaving more general setups for another occasion. On the other hand, in order to render the exposition accessible to readers with various backgrounds (be it probabilistic or ergodic) we made it more expository, and as self-contained as reasonably possible.

Random walks on free groups were the first example that I learned from my teacher Anatoly Moiseevich Vershik when I started to work with him as a second year university student in 1977. He always taught me to look at the \emph{nature of things}. Trying to return to the origins in this article, I dedicate it to Anatoly Moiseevich with immense admiration and gratitude.

\section{Action cocycles and the Radon -- Nikodym problem}

We begin by formulating the general Radon -- Nikodym problem of existence and uniqueness of quasi-invariant probability measures on the boundary $\p\Fc$ with a prescribed Radon~-- Nikodym cocycle. More specifically, we are interested in the \emph{product cocycles}~$r_\pi$ of the boundary action. These are the ones obtained by first endowing the (directed) edges of the Cayley tree of the group $\Fc$ with the weights $\pi_a$ determined by their labels $a\in\Ac$, and then considering the natural regularization of the ratios of infinite products of these weights along the geodesic rays issued from different group elements towards the same boundary point.

\subsection{Action groupoid and its cocycles}

Given a commutative group of values $A$ (for which we use the \emph{additive notation} for the time being), and an action $G\acts X$ of a countable group $G$ on a space~$X$, a function
$$
c:G\times X\to A
$$
is called an ($A$-valued) \textsf{action cocycle} if
\begin{equation} \label{eq:c}
c(g_1 g_2, x) = c(g_2,x) + c(g_1, g_2 x)  \qquad \forall\, g_1,g_2\in G, x\in X \;.
\end{equation}
The \textsf{coboundary} determined by a function $\ph:X\to A$ is a cocycle of the form
$$
c(g,x) = \ph(gx) - \ph(x) \;.
$$

These definitions that go back to \textsc{Mackey} \cite{Mackey66} are best understood in terms of the associated \emph{action groupoid} (also see \textsc{Moore} \cite{Moore08} for a broad contextual and historical discussion). One can think about a \textsf{groupoid} (the most succinctly defined as a \emph{small category with invertible morphisms}) as a collection of arrows \textsf{(morphisms)} between its \textsf{objects}. There is a composition operation that obeys the same rules as the group multipication (groups being precisely the groupoids with only one object). An important difference, though, is that the composition in groupoids is \emph{partial} in the sense that two arrows are composable only if the target of one matches the source of the other one.

\begin{figure}[h]
\begin{center}
\psfrag{x}[][]{$x=g_2^{-1}g_1^{-1}x'$}
\psfrag{y}[][]{$g_2x=g_1^{-1}x'$}
\psfrag{z}[][]{$g_1g_2x=x'$}
\psfrag{a}[cl][cl]{$g_2$}
\psfrag{b}[cl][cl]{$g_1$}
\psfrag{c}[ct][ct]{$g_1g_2$}
\includegraphics[scale=1.]{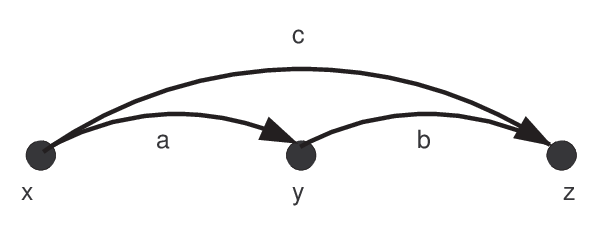}
\end{center}
\caption{Composition in the action groupoid.}
\label{fig:groupoid}
\end{figure}

For the \textsf{action groupoid} $\Gu$ associated with an action $G\acts X$, the space of objects is the action space $X$, and morphisms are the triples $\gu=(x,g,x')$ with $x'=gx$. These morphisms should be thought of as arrows labelled with a group element $g$ and joining a source $x\in X$ with the target $x'=gx$ (the image of $x$ under the action of $g$), see \figref{fig:groupoid}.

Under the identification
$$
\Gu \cong G\times X \;,\qquad (x,g,x')\cong (g,x) \;,
$$
based on the \textsf{source map} $(x,g,x')\mapsto x$, cocycle condition \eqref{eq:c} means that~$c$ is a \textsf{homomorphism} from the action groupoid $\Gu$ to the group $A$:
$$
c(\gu_1 \gu_2) = c(\gu_1) + c(\gu_2)
$$
whenever $\gu_1$ and $\gu_2$ are composable. The coboundary determined by a function $\ph$ on $X$ is then just the difference $\ph(x')-\ph(x)$ of the values of $\ph$ at the target $x'$ and at the source~$x$ of a morphism $(x,g,x')$.

\subsection{Alternative descriptions} \label{sec:alt}

In certain situations (see \secref{sec:rn} below) it is more convenient to use a different parameterization of the action groupoid, the one based on the \textsf{target map} $(x,g,x')\mapsto x'$,
$$
\Gu \cong G\times X \;,\qquad (x,g,x')\cong (g,x') \;.
$$
Under this identification a function $c'$ on $G\times X$ is a \textsf{cocycle} if
\begin{equation} \label{eq:cy}
c'(g_1g_2, x') = c'\!\left(g_2, g_1^{-1} x'\right) + c'(g_1,x')  \qquad \forall\,g_1,g_2\in G, x'\in X \;,
\end{equation}
and a \textsf{coboundary} is a cocycle of the form
\begin{equation} \label{eq:cbdry}
c'(g,x') = \ph(x') - \ph\!\left(g^{-1}x'\right)  = \ph'\!\left(g^{-1} x'\right) - \ph'(x')  \;,
\end{equation}
where the function $\ph'=-\ph$ is the opposite of $\ph$.

$$
\textit{It is definition \eqref{eq:cy} of an action cocycle that we use below.}
$$

\medskip

There is yet another interpretation of action cocycles, in terms of \emph{equivariant families of cocycles} on the group $G$ parameterized by the action space~$X$. Indeed, if $c:G\times X\to A$ is an $A$-valued cocycle as in \dfnref{eq:cy}, then for any $x\in X$
\begin{equation} \label{eq:Cx}
\wt c(x;g_1,g_2) = c(g_1^{-1}g_2,g_1^{-1} x)
\end{equation}
is a \textsf{cocycle} on $G$ in the sense that
\begin{equation} \label{eq:coc}
\wt c(x;g_1,g_2) + \wt c(x;g_2,g_3) = \wt c(x;g_1,g_3) \qquad\forall\,g_1,g_2,g_3\in G \;,
\end{equation}
(in a more formal language, $\wt c(x;\cdot,\cdot)$ is a cocycle of the principal groupoid determined by the equivalence relation whose only class is the whole of $G$), and the family
\begin{equation} \label{eq:wtc}
\wt c=\{\wt c(x;\cdot,\cdot)\}_{x\in X}
\end{equation}
is \textsf{equivariant}, i.e.,
\begin{equation} \label{eq:eq}
\wt c (gx;gg_1,gg_2) = \wt c (x;g_1,g_2) \qquad\forall\,g,g_1,g_2\in G, x\in X \;.
\end{equation}
Conversely, any family $\wt c$ \eqref{eq:wtc} of cocycles on $G$ that satisfies \eqref{eq:eq} gives rise to an action cocycle $c$ by formula
\begin{equation} \label{eq:cwtc}
c(g,x) = \wt c(x;e,g) \;.
\end{equation}

\subsection{Radon -- Nikodym problem} \label{sec:rn}

We remind that a measure $m$ on a space $X$ is called \textsf{quasi-invariant} with respect to a group action $G\acts X$ if all translates $gm, g\in G,$ are equivalent to $m$, or, in other words, if any translate of an $m$-negligible subset is also $m$-negligible.

\begin{dfn}
The \textsf{Radon -- Nikodym cocycle} of a quasi-invariant measure~$m$ on the space $X$ of an
action $G\acts X$ is given by the family of the Radon -- Nikodym derivatives
\begin{equation} \label{eq:De}
\De_m (g,x) = \frac{dgm}{dm}(x) \;, \qquad g\in G, x\in X \;.
\end{equation}
\end{dfn}

The group of values of the Radon -- Nikodym cocycle is the \emph{multiplicative} group $\RR_+^*$ of positive real numbers, and cocycle identity \eqref{eq:cy} is a consequence of the chain rule for the Radon -- Nikodym derivatives. Less formally, the Radon -- Nikodym cocycle is often written as the ``virtual coboundary'', like in \eqref{eq:cbdry}, determined by the infinitesimal ratio
$$
\De_m (g,x) = \frac{dm(g^{-1}x)}{dm(x)} \;,
$$
cf.\ the definition of the Radon -- Nikodym cocycle of a measure quasi-invariant with respect to an equivalence relation due to \textsc{Feldman~--} \linebreak \textsc{Moore}~\cite{Feldman-Moore77}. The corresponding equivariant family \eqref{eq:Cx} of cocycles on~$G$ parameterized by the points from the action space is then
$$
\begin{aligned}
\wt \De_m(x;g_1,g_2) &= \De_m \bigl( g_1^{-1}g_2,g_1^{-1} x \bigr) \\
&= \frac{dg_1^{-1}g_2m}{dm} \bigl( g_1^{-1} x \bigr)
= \frac{dg_2 m}{dg_1 m}(x) \;.
\end{aligned}
$$

The Radon -- Nikodym derivatives \eqref{eq:De} are defined only almost everywhere (``mod~0'') with respect to the measure $m$, and therefore the cocycle identity \eqref{eq:cy} for them is also satisfied only mod 0 for any $g_1, g_2\in G$. However, since the group $G$ is countable, one can choose individual versions of Radon -- Nikodym derivatives and a $G$-invariant full measure subset $X'\subset X$ in such a way that \eqref{eq:cy} be satisfied for \emph{all} $g_1, g_2\in G$ and $x\in X'$.

Given an $\RR^*_+$-valued cocycle $c$ of a group action (more generally, of a groupoid), one can now ask about existence and uniqueness of a measure~$m$ on the action space such that its Radon -- Nikodym cocycle $\De_m$ coincides with $c$. This question is often referred to as the \textsf{Radon -- Nikodym problem} (e.g., see \textsc{Schmidt} \cite[\S 10]{Schmidt77}, or \textsc{Renault} \cite{Renault05}). It plays a central role in the theory of Gibbs measures\,---\,which are essentially defined in terms of their Radon -- Nikodym derivatives (e.g., see \textsc{Haydn -- Ruelle}~\cite{Haydn-Ruelle92}, the discussion of the Anosov -- Sinai cocycle in \textsc{Kaimanovich~-- Lyubich} \cite[Section 2.1.3]{Kaimanovich-Lyubich05}, or the recent comprehensive exposition of \textsc{Paulin -- Pollicott -- Schapira} \cite{Paulin-Pollicott-Schapira15}). Yet another (in fact, very closely related, see \remref{rem:pat}) example is \textsc{Patterson}'s construction \cite{Patterson76} of a measure (``conformal density'') on the limit set of a Fuchsian group (or of more general groups with similar properties) for which the logarithm of the Radon~-- Nikodym cocycle is proportional to the Busemann cocycle, and its various generalizations (cf. \exref{ex:bus} and \exref{ex:unif} below; again see \cite{Paulin-Pollicott-Schapira15} and the numerous references therein for further details).

\section{Free groups: product cocycles of the boundary action \\ and multiplicative Markov measures}

\subsection{Geometry of the Cayley tree} \label{sec:geom}

Let $\Fc=\Fc_k$ be a \textsf{free group} with ${k\ge 2}$ free generators. We denote by $\Ac$ the \textsf{set of generators and their inverses}, and in the usual way identify $\Fc$ with the set of words in the \textsf{alphabet}~$\Ac$ that are \textsf{irreducible} in the sense that any two neighbouring letters are not inverse to each other. By $|g|$ we denote the \textsf{length} of the irreducible word representing a group element $g\in \Fc$.

The Cayley graph of $\Fc$ with respect to the generating set $\Ac$ is a homogeneous tree of degree $2k$, and the \textsf{graph distance} $d(g_1,g_2)$ between two vertices~$g_1,g_2\in\Fc$ is the length of~$g_1^{-1} g_2$. By $[g_1,g_2]$ we denote the unique \textsf{geodesic} joining $g_1$ and $g_2$. The elements of the alphabet $\Ac$ provide a labelling of the directed edges of the Cayley graph, and the sequence of these labels along the geodesic $[g_1,g_2]$ is precisely the word $g_1^{-1}g_2$.

The group $\Fc$ has a natural \textsf{compactification}
$$
\ol \Fc = \Fc \cup \p \Fc
$$
which admits a number of equivalent descriptions (e.g., see \textsc{Kapovich -- Benakli} \cite{Kapovich-Benakli02}). In particular, the \textsf{boundary} $\p\Fc$ can be identified with the set of the right \textsf{infinite irreducible words} in the alphabet $\Ac$, and in these terms convergence of a sequence $g_n\in\Fc$ to $\ga\in\p\Fc$ just means that the words~$g_n$ converge to $\ga$ pointwise (i.e., for any $m>0$ the initial length~$m$ segments \textsf{(truncations)}~$[g_n]_m$ of the words~$g_n$ converge to the $m$-truncation $[\ga]_m$ of $\ga$). The action of the group~$\Fc$ on itself extends to a continuous left \textsf{boundary action}
$$
\Fc\acts\p\Fc: (g,\ga)\mapsto g\ga \;,
$$
where $g\ga$ is the infinite word obtained by the concatenation of $g$ and $\ga$ followed, if necessary, by cancellations (in the same way as for the group operation in $\Fc$).

\begin{figure}[h]
\begin{center}
\psfrag{e}[cr][cr]{$e$}
\psfrag{a}[cl][cl]{$\p\Fc$}
\psfrag{o}[cl][cl]{$\ga$}
\includegraphics[scale=.6]{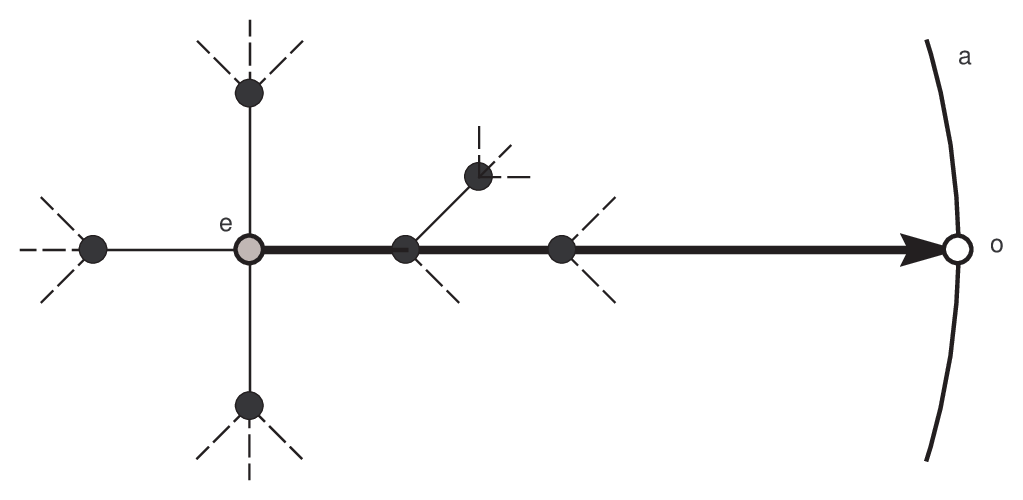}
\end{center}
\caption{The Cayley tree and its boundary.}
\label{fig:graph}
\end{figure}

Any geodesic ray of the graph metric converges in the above compactification, and for any starting point $g\in\Fc$ and a boundary point $\ga\in\Fc$ there is a unique ray $[g,\ga)$ convergent to $\ga$. Reading the edge labels along the ray $[e,\ga)$ issued from the group identity produces then the infinite irreducible word $\ga$, see \figref{fig:graph}.

\pagebreak

\begin{figure}[h]
\begin{center}
\psfrag{a}[cl][cl]{$\p\Fc$}
\psfrag{o}[cl][cl]{$\ga$}
\psfrag{x}[cr][cr]{$g_1$}
\psfrag{y}[cr][cr]{$g_2$}
\psfrag{z}[cr][cr]{$g_1 \wedge_\ga g_2$}
\includegraphics[scale=.6]{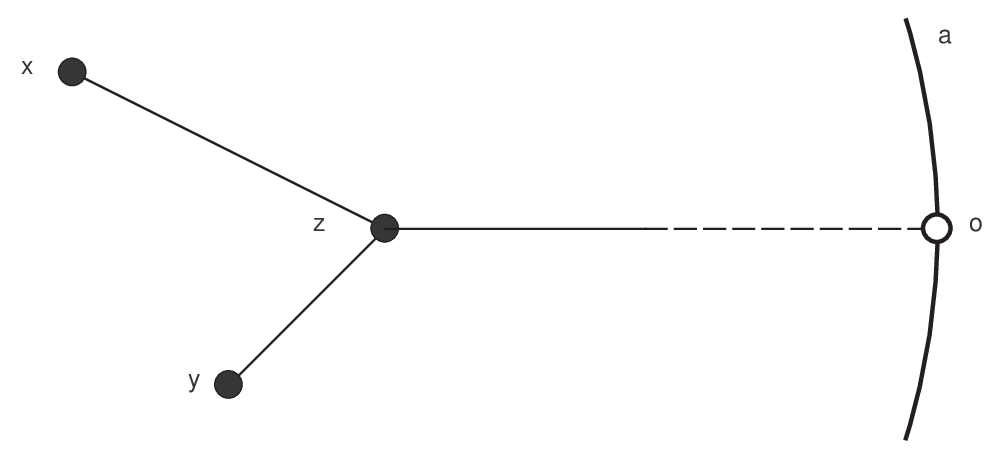}
\end{center}
\caption{The confluent of two group elements with respect to a boundary point.}
\label{fig:bus}
\end{figure}

Any two geodesic rays $[g_1,\ga), [g_2,\ga)$ with a common boundary point $\ga$ coincide starting from a point called their \textsf{confluent} and denoted $g_1 \wedge_\ga g_2$, see \figref{fig:bus}.

For any $g\in \Fc$ we denote by $\Sc_g\subset \p \Fc$ the \textsf{cylinder set} of all infinite irreducible words that begin with $g$ (in particular, $\Sc_e=\p\Fc$ for the group identity $e$). Geometrically, $\Sc_g$ is the set of all $\ga\in\p\Fc$ such that the geodesic ray $[e,\ga)$ passes through $g$, i.e., the \textsf{shadow} of $g$ with respect to a light source at the group identity.

\begin{figure}[h]
\begin{center}
\psfrag{e}[][]{$g$}
\psfrag{g}[][]{$ga$}
\psfrag{a}[cl][cl]{$g\Sc_a$}
\psfrag{b}[cr][cr]{$\p\Fc\sm g \Sc_a=(ga) \Sc_{a^{-1}}$}
\hskip 2cm\includegraphics[scale=.6]{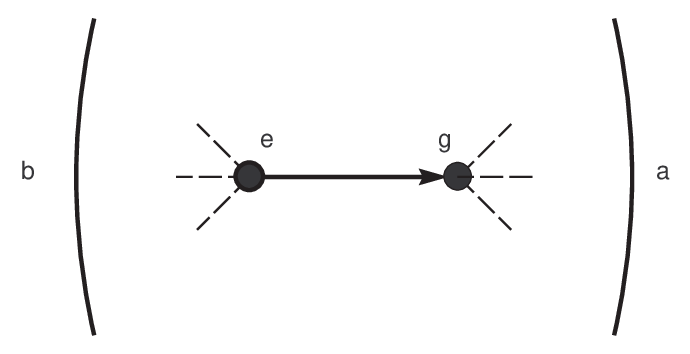}
\end{center}
\caption{The partition of the boundary $\p\Fc$ determined by a directed edge of the Cayley graph.}
\label{fig:two}
\end{figure}

Removing a directed edge $[g,ga]$ (where $g\in\Fc$ and $a\in\Ac$) splits the Cayley tree of the group $\Fc$ into two connected components. It also gives rise to a partition of the boundary $\p\Fc$ into two subsets (the boundaries of these components): the ``positive'' one~$g\Sc_a$ is the set of limit points of the geodesic rays issued from the initial point $g$ of the edge and passing through its endpoint~$ga$, whereas, \emph{mutatis mutandis}, the ``negative'' one $\p\Fc\sm g\Sc_a = (ga) \Sc_{a^{-1}}$ is the set of limit points of the geodesic rays issued from the endpoint $ga$ of the edge and passing through the initial point $g$, see \figref{fig:two}. Accordingly, we say that the orientation of a directed edge $[g,ga]$ is \textsf{positive} (resp., \textsf{negative}) with respect to a boundary point $\ga$ if $\ga\in g\Sc_a$ (resp., if $\ga\in \Fc\sm g\Sc_a = (ga) \Sc_{a^{-1}}$).

\subsection{The potential of a boundary action cocycle}

\begin{dfn} \label{dfn:potential}
The \textsf{potential} of an $A$-valued cocycle~$c$ of the boundary action $\Fc\acts\p\Fc$, or, equivalently, of the associated equivariant family \eqref{eq:Cx} of cocycles $\{\wt c(\ga;\cdot,\cdot)\}_{\ga\in\p\Fc}$ on $\Fc$, is the function
\begin{equation} \label{eq:xi}
\upxi:\p\Fc\to A\;, \qquad \upxi(\ga) = c\bigl([\ga]_1,\ga\bigr) = \wt c \bigl(\ga;e,[\ga]_1\bigr) \;,
\end{equation}
where $[\ga]_1$ denotes the first letter of an infinite word $\ga\in\p\Fc$. In more geometrical terms $\upxi(\ga)$ is the value of the cocycle $\wt c(\ga,\cdot,\cdot)$ on the initial edge $\bigl[e,[\ga]_1\bigr]$ of the geodesic ray~$[e,\ga)$.
\end{dfn}

\begin{prp}
Any function $\upxi:\p\Fc\to A$ uniquely extends to an $A$-valued cocycle of the boundary action $\Fc\acts\p\Fc$ with the potential $\upxi$.
\end{prp}

\begin{proof}
It is more convenient to argue first in terms of the $\Fc$-equivariant families of cocycles~$\wt c$ on $\Fc$ parameterized by $\p\Fc$ (see \secref{sec:alt}) instead of action cocycles $c$. The Cayley graph of the group~$\Fc$ being a tree, it has no loops, and therefore in order to define a cocycle on $\Fc$ it is enough to specify just its values on directed edges without imposing any additional conditions other than requiring that the cocycle take opposite values for two opposite orientations of each edge.

By the definition of a potential, we have to put
$$
\wt c (\ga; e, a) = \upxi (\ga) \qquad\forall\,a\in\Ac, \ga\in\Sc_a \;,
$$
and by equivariance then
$$
\wt c (\ga; g, ga) = \upxi (g^{-1}\ga)
$$
for all triples $(\ga;g,ga)$ such that the orientation of the edge $[g,ga]$ is positive with respect to $\ga$, i.e., $\ga\in g\Sc_a$. Finally, if the orientation of $[g,ga]$ with respect to $\ga$ is negative, i.e., if $\ga\in\p\Fc\sm g \Sc_a=(ga) \Sc_{a^{-1}}$ (see \figref{fig:two}), then one has to reverse the edge and put
$$
\wt c (\ga; g, ga) = - \wt c (\ga; ga, g) = - \wt c \bigl(\ga; ga, (ga) a^{-1}\bigr) = - \upxi \bigl( (ga)^{-1}\ga \bigr) \;.
$$
Thus, $\upxi$ uniquely determines an equivariant family $\wt c$ of cocycles on $\Fc$.

In view of formula \eqref{eq:cwtc}, the corresponding boundary action cocycle $c$ with the potential~$\upxi$ is defined then by putting for any $a\in\Ac$
$$
c(a,\ga) =
\begin{cases}
\upxi(\ga) \;, & \ga \in \Sc_a \;, \\
-\upxi(a^{-1}\ga) \;, & \ga \in \p\Fc\sm\Sc_a \;,
\end{cases}
$$
and further using the chain rule \eqref{eq:cy} for defining $c(g,\ga)$ with $|g|>1$.
\end{proof}

\begin{rem}
In terms of the formal series
\begin{equation} \label{eq:Xi}
\Xi(\ga) = \mathquote{\sum\limits_{n=0}^\infty \upxi(T^n\ga)} \;,
\end{equation}
where
\begin{equation} \label{eq:T}
T: \p\Fc\to\p\Fc\;,\qquad  \ga \mapsto [\ga]_1^{-1}\ga \;,
\end{equation}
is the \textsf{shift} on $\p\Fc$ (which consists in erasing the first letter of the corresponding infinite word), the cocycles $\wt c(\ga,\cdot,\cdot)$ are regularizations of the formal expressions
\begin{equation} \label{eq:Cga}
\wt c (\ga;g_1,g_2) = \mathquote{\,\Xi(g_1^{-1}\ga) - \Xi(g_2^{-1}\ga)\,} \;.
\end{equation}
Indeed, since the forward $T$-orbits of $g_1^{-1}\ga$ and $g_2^{-1}\ga$ eventually merge, the corresponding sums \eqref{eq:Xi} differ only by finitely many summands, see \figref{fig:bus}.
\end{rem}

\begin{ex} \label{ex:bus}
The $\ZZ$-valued cocycle $b$ of the boundary action determined by the constant potential
\begin{equation} \label{eq:be}
\upbeta(\ga) = -1
\end{equation}
is called the \textsf{Busemann cocycle}. Geometrically, it corresponds to the equivariant family of $\ZZ$-valued Busemann cocycles on $\Fc$
$$
\wt b (\ga;g_1, g_2) = d(g_2, g_1 \wedge_\ga g_2) - d(g_1, g_1 \wedge_\ga g_2) \;, \qquad \ga\in\p\Fc \;,
$$
defined by regularizing the formal expression
$$
\wt b (\ga;g_1,g_2) = \mathquote{\,d(g_2,\ga)-d(g_1,\ga)\,}
$$
(see \figref{fig:bus}, cf. \textsc{Kaimanovich} \cite{Kaimanovich90} and \textsc{Cartwright -- Kai\-ma\-no\-vich~-- Woess} \cite{Cartwright-Kaimanovich-Woess94}), which means that the cocycle $\wt b(\ga;\cdot,\cdot)$ can be considered as the discrete gradient of the ``distance'' to $\ga\in\p\Fc$ (which is why the sign in \eqref{eq:be} is negative). The \textsf{kernel} of the cocycle $\wt b(\ga;\cdot,\cdot)$
$$
\bigl\{ (g_1,g_2)\in \Fc\times\Fc: \wt b(\ga;g_1,g_2)=0 \bigr\}
$$
is an equivalence relation on $\Fc$, the classes of which (i.e., the subsets of~$\Fc$ that consist of all points ``equidistanced'' from $\ga$) are the \textsf{horospheres} centered at $\ga$.
\end{ex}

\begin{ex} \label{ex:unif}
The \textsf{uniform} probability measure $\me$ on $\p\Fc$ is the one defined by requiring that the measures of the shadows of all group elements of the same length be equal, i.e., that
$$
\me(\Sc_g) = \frac1{2k (2k-1)^{|g|-1}} \qquad \forall\, g\in\Fc\sm\{e\} \;,
$$
where $k$ is the number of generators of $\Fc$. It is well-known that the measure~$\me$ is quasi-invariant with respect to the boundary action, and its Radon~-- Nikodym cocycle is
$$
\De_{\textstyle\matheu{m}} (g,\ga) = (2k-1)^{-b(g,\ga)} \qquad\forall\,g\in\Fc,\ga\in\p\Fc \;.
$$
Indeed, it is enough to verify that the potential of the Radon -- Nikodym cocycle is the constant $2k-1$, i.e., to check that $\De_{\textstyle\matheu{m}} (a,\ga)=2k-1$ whenever $a\in\Ac$ and $\ga\in\Sc_a$. The measure $a\me$ is uniform for the observer looking at the boundary $\p\Fc$ from the point $a$, whence a comparison of the shadow measures yields the claim, see \figref{fig:two}.
\end{ex}

\subsection{The Radon -- Nikodym problem for product cocycles} \label{sec:RN}

\begin{dfn} \label{dfn:product}
Given an assignment
\begin{equation} \label{eq:pi}
\pi=\{\pi_a\}_{a\in\Ac}
\end{equation}
of positive \textsf{weights} to the generators and their inverses, the associated $\RR^*_+$-valued \textsf{product cocycle} $r=r_\pi$ of the boundary action $\Fc\acts\p\Fc$ is the one determined by the potential
$$
\uprho(\ga) = \frac1{\pi_{[\ga_1]}} \;
$$
i.e., for any $a\in\Ac$
\begin{equation} \label{eq:r}
r(a,\ga) =
\begin{cases}
\dfrac1{\pi_a} \;, & \ga \in \Sc_a \;, \\ \\
\pi_{a^{-1}} \;, & \ga \in \p\Fc\sm\Sc_a
\end{cases}
\end{equation}
In the same way as above, by
$$
\wt r(\ga;g_1,g_2)= r\!\left(g_1^{-1}g_2,g_1^{-1}\ga\right)
$$
we denote the corresponding equivariant family of cocycles on $\Fc$ indexed by $\ga\in\p\Fc$.
\end{dfn}

Assignment \eqref{eq:pi} allows one to label the directed edges of the Cayley graph with the corresponding weights. Given two elements $g_1,g_2\in\Fc$, we denote by $\pi[g_1, g_2]$ the product of these weights along the geodesic $[g_1,g_2]$. Then
$$
\frac{\pi[g_2,g]}{\pi[g_1,g]} = \frac{\pi[g_2,g_1 \wedge_g g_2 ]}{\pi[g_1,g_1 \wedge_g g_2 ]}
\qquad\forall\,g,g_1,g_2\in G \;,
$$
where $g_1 \wedge_g g_2$ denotes the confluent (the first common point) of the geodesic segments $[g_1,g]$ and $[g_2,g]$ (cf.\ \figref{fig:bus}), and
$$
\wt r(\ga;g_1,g_2) = \frac{\pi[g_2,g_1 \wedge_\ga g_2 ]}{\pi[g_1,g_1 \wedge_\ga g_2 ]}
$$
is a regularization of the expression
$$
\wt r(\ga;g_1,g_2) = \text{\raisebox{7pt}{``}}\, \frac{\pi[g_2,\ga)}{\pi[g_1,\ga)}\, \text{\raisebox{7pt}{''}} \;,
$$
where $\pi[g,\ga)$ is the formal product of the weights of the edges of the geodesic ray $[g,\ga)$, cf.~\eqref{eq:Cga}. In particular, if all weights $\{\pi_a\}_{a\in\Ac}$ have the same value~$p>0$, then the associated product cocycle is $r=p^b$, where~$b$ is the Busemann cocycle of the boundary action, cf. \exref{ex:unif}.

\begin{thm} \label{thm:final}
Given a collection of positive weights $\pi=\{\pi_a\}_{a\in \Ac}$ on the alphabet $\Ac$, the associated product cocycle $r=r_\pi$ of the boundary action $\Fc\acts\p\Fc$ can be realized as the Radon~-- Nikodym cocycle $\De_\nub$ of a quasi-invariant probability measure $\nub$ on $\p \Fc$ if and only if
\begin{equation} \label{eq:pia1}
\pi_a<1 \qquad\forall\,a\in\Ac \;,
\end{equation}
and
\begin{equation} \label{eq:sum}
\sum_{a\in\Ac} \frac{\pi_a(1-\pi_{a^{-1}})}{1-\pi_a\pi_{a^{-1}}} = 1 \;.
\end{equation}
Under these conditions the measure $\nub$ is unique, and it is determined by its values on the shadow sets
\begin{equation} \label{eq:nus}
\nub(\Sc_g) = \pi_{a_1}\dots\pi_{a_{n-1}}\nu_{a_n} \qquad \forall\, g=a_1\dots a_n \in \Fc \;,
\end{equation}
where
\begin{equation} \label{eq:nuaa}
\nu_a = \nub(\Sc_a) = \frac{\pi_a(1-\pi_{a^{-1}})}{1-\pi_a\pi_{a^{-1}}} \;, \qquad a\in\Ac \;,
\end{equation}
are the weights of the probability measure $\nu=\{\nu_a\}$ describing the first letter of $\nub$-distributed infinite irreducible words.
\end{thm}

\begin{proof}
Let $\nub$ be a quasi-invariant probability measure on $\p \Fc$ whose Radon~-- Nikodym cocycle is $\De_\nub=r$. Since the boundary action of $\Fc$ is minimal, the measure $\nub$ is fully supported, so that the measures $\nu_a=\nub(\Sc_a)$ of the shadows $\Sc_a$ are all strictly positive, and obviously
\begin{equation} \label{eq:sums}
\sum_{a\in\Ac} \nu_a = 1 \;.
\end{equation}
Then by \eqref{eq:r} for any $a\in\Ac$
\begin{equation} \label{eq:nua}
\begin{aligned}
\nu_a = \nub(\Sc_a) &= \pi_a \cdot a\nub(\Sc_a) = \pi_a \cdot \nub(a^{-1}\Sc_a) \\
&= \pi_a \cdot \nub (\p\Fc\sm\Sc_{a^{-1}}) \\
&= \pi_a (1-\nu_{a^{-1}}) \;,
\end{aligned}
\end{equation}
so that
\begin{equation} \label{eq:pia}
\pi_a = \frac{\nu_a}{1-\nu_{a^{-1}}} \;,
\end{equation}
and therefore condition \eqref{eq:pia1} is satisfied. After substituting the analogous expression for~$\nu_{a^{-1}}$ in terms of $\nu_a$, \eqref{eq:nua} yields
\begin{equation} \label{eq:piaex}
\nu_a = \pi_a(1-\pi_{a^{-1}}(1-\nu_a)) \;,
\end{equation}
whence \eqref{eq:nuaa}, and \eqref{eq:sums} implies that condition \eqref{eq:sum} is satisfied.

Then for any $g=a_1\dots a_n\in \Fc$ with $|g|=n>1$ one obtains \eqref{eq:nus} in the same way as in~\eqref{eq:nua}:
$$
\begin{aligned}
\nub(\Sc_{a_1\dots a_n}) &= \pi_{a_1} \cdot a_1\nub(\Sc_{a_1\dots a_n})
= \pi_{a_1} \cdot \nub(a_1^{-1}\Sc_{a_1\dots a_n})
= \pi_{a_1} \cdot \nub(\Sc_{a_2\dots a_n}) \\
&= \pi_{a_1}\ldots\pi_{a_{n-1}}\nub(\Sc_{a_n})
= \pi_{a_1}\ldots\pi_{a_{n-1}}\nu_{a_n} \;,
\end{aligned}
$$
whence $\nub$ is, indeed, the unique measure with the Radon -- Nikodym cocycle $\De_\nub=r$.

Conversely, given a collection of weights satisfying \eqref{eq:pia1} and \eqref{eq:sum}, the values of the measures of cylinder sets produced by formulas \eqref{eq:nus} and \eqref{eq:nuaa} are Kolmogorov consistent, because for any $a\in \Ac$
$$
\begin{aligned}
\sum_{b\in\Ac:ab\neq e} \pi_a\nu_b
&= \pi_a(1-\nu_{a^{-1}}) \\
&= \pi_a \left( 1 - \frac{\pi_{a^{-1}}(1-\pi_a)}{1-\pi_a\pi_{a^{-1}}} \right)
= \frac{\pi_a(1-\pi_{a^{-1}})}{1-\pi_a\pi_{a^{-1}}} = \nu_a
\;,
\end{aligned}
$$
and therefore they uniquely determine a probability measure $\nub$ on $\p \Fc$. This measure is then clearly quasi-invariant with the Radon -- Nikodym cocycle $\De_\nub=r$.
\end{proof}

\subsection{Quasi-product Markov measures}

The irreducibility condition \linebreak evoked in \secref{sec:geom} (mutually inverse generators are not allowed to stand next to each other) is an example of a structure well-known in symbolic dynamics.

A \textsf{topological Markov chain} $\Ac_\Upsi^\NN$ over a finite state space $\Ac$ is determined by an \mbox{$\Ac\times\Ac$} \textsf{admissibility matrix} $\Upsi$ with $\{0,1\}$ entries: it is the space of all sequences $(a_1,a_2,\dots)$ with entries $a_t\in\Ac$ subject to the condition
$$
\Upsi_{a_t a_{t+1}}=1 \qquad\forall\,t\ge 1 \;,
$$
and endowed with the compact topology of coordinate-wise convergence.\footnotemark\ Given an initial distribution $\nu$ on $\Ac$ and a stochastic \textsf{transition matrix (Markov operator)}~$P$ over $\Ac$ \textsf{subordinated} to $\Upsi$ (i.e., such that $P_{ab}=0$ whenever $\Upsi_{ab}=0$) whose entries are interpreted as \textsf{transition probabilities} $p_{a\lt b} = P_{ab}$, the associated \textsf{Markov measure} $\nub$ on $\Ac_\Upsi^\NN$ is defined in the usual way by its values on the cylinder sets
\begin{equation} \label{eq:M}
\nub \left( C_{a_1 \dots a_t} \right) = \nu_{a_1} p_{a_1\lt a_2} \dots p_{a_{t-1}\lt a_t} \;.
\end{equation}

\footnotetext{\;Usually the time $t$ in a topological Markov chain runs from $t_0=0$. However, since our main example is the boundary of the free group, in this situation it is more natural to enumerate the letters of infinite irreducible words from $\p\Fc$ starting from $t_0=1$.}

However, by using the admissibility matrix $\Upsi$, one can also produce a transition matrix $P=P_{\si,\Upsi}$ just from a single distribution $\si=(\si_a)$ on $\Ac$ by putting
\begin{equation} \label{eq:pab}
p_{a\lt b} = \frac{\si_b}{\si(\Upsi_{a\lt})} = \frac{\si_b}{(\Upsi\si)_a} \;,
\end{equation}
where
$$
\Upsi_{a\lt} = \{b\in\Ac:\Upsi_{ab}=1\}
$$
is the set of all states to which the topological Markov chain is allowed to transition from~$a$, and $\Upsi\si$ is the result of the action of the matrix $\Upsi$ on the vector $\si$ (considered here as a function, so that it is treated as a column vector and stands to the right of the matrix~$\Upsi$). In other words, for any $a\in\Ac$ the associated transition distribution is the normalized restriction of the measure $\si$ to the set $\Upsi_{a\lt}$, and formula \eqref{eq:M} for the resulting Markov measure $\nub$ with the initial distribution $\nu$ and the transition probabilities \eqref{eq:pab} takes the form
\begin{equation} \label{eq:nC}
\begin{aligned}
\nub \left( C_{a_1 \cdot \ldots \cdot a_t} \right)
&= \nu_{a_1} \cdot \frac{\si_{a_2} \cdot \ldots \cdot \si_{a_t}}{(\Upsi\si)_{a_1} \cdot \ldots \cdot (\Upsi\si)_{a_{t-1}}} \\
&= \frac{\nu_{a_1}}{\si_{a_1}} \cdot \frac{\si_{a_1} \cdot \ldots \cdot \si_{a_t}}{(\Upsi\si)_{a_1} \cdot \ldots \cdot (\Upsi\si)_{a_t}} \cdot (\Upsi\si)_{a_t} \\
&= \al_{a_1} \cdot \rho_{a_1} \cdot \ldots \cdot \rho_{a_t} \cdot \om_{a_t} \;,
\end{aligned}
\end{equation}
where $a_1,a_2,\dots,a_t$ is a $\Upsi$-admissible sequence of length $t\ge 2$, and
\begin{equation} \label{eq:rho}
\rho_a = \frac{\si_a}{(\Upsi\si)_a} \;, \quad \al_a = \frac{\nu_a}{\si_a} \;, \quad \om_a = (\Upsi\si)_a  \qquad \forall\;a\in\Ac \;.
\end{equation}

In other words, the finite dimensional distributions of $\nub$ differ from the restrictions of the respective product measures $\rho^{\otimes t}$ to the corresponding sets of length $t$ admissible finite paths just by the correction densities $\al$ and~$\om$ that only depend on $a_1$ and $a_t$, respectively. Since the transition probabilities~\eqref{eq:pab} and the weights \eqref{eq:rho} do not change when the measure~$\si$ is multiplied by a constant, without loss of generality we may assume that $\si$ is normalized and its total mass is 1. It is worth noticing that the measure~$\rho$~\eqref{eq:nC} on the alphabet $\Ac$ is not a probability one (except for very degenerate exceptional cases), and therefore the infinite product of $\rho$ by itself is not well-defined.

\begin{dfn} \label{dfn:qp}
A Markov measure $\nub$ on a topological Markov chain $\Ac_\Upsi^\NN$ is called \textsf{quasi-product} if its transition probabilities are determined by a probability \textsf{transition measure} $\si$ and the admissibility matrix $\Upsi$ by formula~\eqref{eq:pab}, or, equivalently, if its values on cylinder sets are given by formula~\eqref{eq:nC}, where $\nu$ is the initial distribution of $\nub$. If, moreover, the initial and the transition distributions coincide, then we say that $\nub$ is a \textsf{multiplicative Markov} measure with the \textsf{base distribution} $\nu=\si$.
\end{dfn}

In our situation $\Ac$ is the set of free generators and their inverses of the free group $\Fc$, and the admissibility matrix is
\begin{equation} \label{eq:SiF}
\Upsi_{ab} =
\begin{cases}
1 \;,\quad ab\neq e \;, \\
0 \;,\quad ab = e \;,
\end{cases}
\end{equation}
so that the arising topological Markov chain $\Ac_\Upsi^\NN$ can be identified with the boundary $\p\Fc$ (the set of right infinite irreducible words in the alphabet $\Ac$). As a dynamical system, the topological Markov chain $\Ac_\Upsi^\NN\cong\p\Fc$ is endowed with the action of the shift $T$ \eqref{eq:T}.

Then
$$
\Upsi_{a\lt} = \Ac \,\setminus\, \{a^{-1}\} \;,
$$
and
$$
(\Upsi\si)_a = \sum_{b\in\Ac} \si_b - \si_{a^{-1}} = 1 - \si_{a^{-1}} \;,
$$
so that the transition probabilities \eqref{eq:pab} take the form
$$
p_{a\lt b} = \frac{\si_b}{1-\si_{a^{-1}}} \;.
$$

Returning to the measure $\nub$ on $\p\Fc$ from \thmref{thm:final}, we now see that it is Markov with the initial distribution $\{\nu_a\}_{a\in\Ac}$ and the transition probabilities
\begin{equation} \label{eq:P}
p_{a\lt b} = \frac{\pi_a\nu_b}{\nu_a} \;,
\end{equation}
subordinated to the admissibility matrix \eqref{eq:SiF}. Moreover, since the above formula for the transition probabilities amounts to dividing $\nu_b$ by the normalizing constant $\nu_a/\pi_a$ (which is then automatically equal to $1 - \nu_{a^{-1}}$), the measure $\nub$ is actually a multiplicative Markov one determined by the base distribution $\nu=\{\nu_a\}_{a\in\Ac}$ and the admissibility matrix $\Upsi$ \eqref{eq:SiF}. Thus, \thmref{thm:final} can be reformulated in the following way:

\begin{thm} \label{thm:qi}
The class of quasi-invariant probability measures on $\p\Fc$ with a product Radon -- Nikodym cocycle coincides with the class of multiplicative Markov measures. More precisely:

\begin{enumerate}[label = \upshape(\roman*)]
\item
Given a collection of weights $\pi=\{\pi_a\}_{a\in \Ac}$, the associated product cocycle $r_\pi$ (\dfnref{dfn:product}) can be realized as the Radon -- Nikodym cocycle of a quasi-invariant probability measure $\nub$ on $\p \Fc$ if and only if conditions \eqref{eq:pia1} and \eqref{eq:sum} are satisfied. In this case the measure~$\nub$ is unique, and it is the multiplicative Markov measure with the base distribution $\nu=\{\nu_a\}_{a\in \Ac}$ \eqref{eq:nuaa}.
\item
Conversely, given a probability multiplicative Markov measure $\nub$ \linebreak on~$\p \Fc$ with the base distribution $\nu=\{\nu_a\}$, its Radon -- Nikodym \linebreak cocycle~$\De_\nub$ is the product cocycle~$r_\pi$ determined by the weights \linebreak ${\pi=\{\pi_a\}_{a\in \Ac}}$ \eqref{eq:pia} that satisfy conditions \eqref{eq:pia1} and~\eqref{eq:sum}.
\end{enumerate}
\end{thm}

\subsection{Perron -- Frobenius interpretation}

Condition \eqref{eq:sum} on the \linebreak weights $\{\pi_a\}$ that appears in \thmref{thm:final} and \thmref{thm:qi} admits a very natural interpretation in terms of the Perron -- Frobenius theorem. We denote by $\rb(M)$ the \textsf{Perron -- Frobenius eigenvalue} ($\equiv$ the spectral radius) of an irreducible square matrix $M$ with non-negative entries (e.g., see \textsc{Seneta}~\cite{Seneta06}).

\begin{thm} \label{thm:pf}
Let $\{\pi_a\}_{a\in\Ac}$ be a collection of positive weights on the alphabet $\Ac$, and let~$\Uppi$ be the corresponding diagonal matrix. Then $\rb(\Uppi\Upsi)=1$ if and only if conditions \eqref{eq:pia1} and \eqref{eq:sum} are satisfied.
\end{thm}

\begin{proof}
In one direction this is a direct consequence of \thmref{thm:qi}(i). Indeed, assuming that conditions \eqref{eq:pia1} and \eqref{eq:sum} are satisfied, let $P$ be the arising stochastic matrix ($\equiv$ Markov operator) whose entries are the transition probabilities \eqref{eq:P}. In the matrix form
\begin{equation} \label{eq:coh}
P = \Upnu^{-1} \Uppi \Upsi \Upnu \;,
\end{equation}
where $\Upnu$ is the diagonal matrix with the entries $\{\nu_a\}_{a\in\Ac}$ (the appearance of~$\Upsi$ in formula \eqref{eq:coh} is due to the fact that transitions \eqref{eq:P} are subject to the admissibility condition $ab\neq e$ represented by $\Upsi$), i.e., the operator $P$ is \emph{cohomologous} to the operator $\Uppi\Upsi$.\footnotemark\ Since $\rb(P)=1$ (because $P$ is Markov), we conclude that $\rb(\Uppi\Upsi)=1$.

\footnotetext{\; The transformations of the form \eqref{eq:coh} are well-known in the Markov theory under the name of \emph{$h$-transforms,} where~$h$ is the standard notation for the function on the state space that establishes the cohomology (and which is required to be harmonic for the original Markov operator to guarantee that the resulting operator be Markov as well). This construction was first introduced by \textsc{Feller} \cite{Feller56}, and then prominently used by \textsc{Doob} for describing conditional Markov processes \cite{Doob57}, which is why it is often called \emph{Doob's $h$-transform}, or just \emph{Doob's transform} (although Doob himself explicitly acknowledged the work of Feller and the influence of Brelot's relative potential theory in \cite[pp.~441--442]{Doob59}). In our situation the role of $h$ is played by the family $\{\nu_a\}_{a\in\Ac}$. However, in the context of \thmref{thm:pf} it is only the operator~$P$ that is Markov, whereas the operator $\Uppi\Upsi$ is not.}

Conversely, let us first note that if $\pi_a\ge 1$ for a certain $a\in\Ac$, then for any $b\neq a,a^{-1}$ already
$$
\rb
\begin{pmatrix}
\pi_a & \pi_a \\
\pi_b & \pi_b
\end{pmatrix} > 1
$$
for the $2\times 2$ submatrix of $\Uppi\Upsi$ over the states $a$ and $b$, and therefore ${\rb(\Uppi\Upsi)>1}$. Thus, if $\rb(\Uppi\Upsi)=1$, then condition \eqref{eq:pia1} is satisfied. We can now explicitly solve the equation
$$
\ph\Uppi\Upsi=\ph
$$
for the \textsf{left Perron -- Frobenius eigenvector} $\ph=\{\ph_a\}_{a\in\Ac}$ of the matrix $\Uppi\Upsi$ (which is unique up to a multiplicative constant). Namely,
$$
\ph_b = \sum_{a:a\neq b^{-1}} \pi_a\ph_a = S - \pi_{b^{-1}} \ph_{b^{-1}} \qquad \forall\, b\in\Ac \;,
$$
where
$$
S = \sum_{a\in\Ac} \pi_a \ph_a \;.
$$
The latter system splits into pairs of equations corresponding to mutually inverse generators, cf. \eqref{eq:piaex}, whence
\begin{equation} \label{eq:ph}
\ph_b = \frac{1-\pi_{b^{-1}}}{1-\pi_b\pi_{b^{-1}}} S \qquad\forall\,b\in\Ac \;,
\end{equation}
which in view of the definition of the sum $S$ leads to condition \eqref{eq:sum}.
\end{proof}

\begin{rem}
The Perron -- Frobenius eigenvalues of the matrices $\Uppi\Upsi$ and~$\Upsi\Uppi$ coincide. There are several reasons for that, the most general one being that the Perron -- Frobenius eigenvalue of the product of two matrices does not depend on the order in which they are multiplied. Moreover, since $\Uppi$ is invertible, $\Uppi\Upsi$ and $\Upsi\Uppi$ are cohomologous. Yet another reason is that in our situation the matrix $\Upsi$ is symmetric, so that $\Upsi\Uppi$ and $\Uppi\Upsi$ are mutually transpose. The latter observation implies that if $\ph$ is a left Perron -- Frobenius eigenvector of $\Uppi\Upsi$ \eqref{eq:ph}, then $\Upsi\Uppi\ph=\ph$, whence $\Uppi\Upsi\Uppi\ph=\Uppi\ph$, i.e., $\Uppi\ph$ is a \emph{right} Perron -- Frobenius eigenvector of $\Uppi\Upsi$.
\end{rem}

\begin{rem} \label{rem:PF}
The spectral characterization of the product cocycles for which the Radon -- Nikodym problem is solvable that results from the combination of \thmref{thm:qi} and \thmref{thm:pf} can be considered as a particular case of the general \emph{Ruelle -- Perron -- Frobenius theorem} for H\"older potentials on topological Markov chains (e.g., see \textsc{Parry -- Pollicott} \cite[Theorem 2.2 and Remark 1 from Chapter 2]{Parry-Pollicott90}).
\end{rem}

\begin{rem} \label{rem:pat}
There is a general approach to the Radon -- Nikodym \linebreak problem for boundary actions of hyperbolic spaces due to \textsc{Patterson} \cite{Patterson76} (also see \textsc{Sullivan} \cite{Sullivan79}, \textsc{Kaimanovich -- Lyubich} \cite{Kaimanovich-Lyubich05}, or \textsc{Denker~--} \linebreak \textsc{Stratmann}~\cite{Denker-Stratmann12} for further generalizations and developments). In our context it takes the following form. Given a collection of positive weights $\{\pi_a\}_{a\in\Ac}$ and a number $\de\ge 0$, let us consider a measure (\emph{a priori} not necessarily finite) $\eta_\de$ on $\Fc$ with the weights
$$
\eta_\de(g) = \pi_g^\de \;, \qquad g\in\Fc \;,
$$
where $\pi_g = \pi[e,g]$ is the product of the weights $\pi_a$ along the edges of the geodesic $[e,g]$ that joins the group identity with $g$. Then there exists a critical value $\de_0$ such that $\|\eta_\de\|=\infty$ for $0\le\de\le\de_0$, and $\|\eta_\de\|<\infty$ for $\de>\de_0$. The normalized measures $\eta_\de/\|\eta_\de\|$ weakly converge, as $\de$ tends to $\de_0$ from above, to a probability measure $\eta$ on $\p\Fc$ whose Radon -- Nikodym cocycle is precisely the product cocycle corresponding to the weights $\{\pi_a^{\de_0}\}_{a\in\Ac}$. Thus, $\de_0=1$ if and only if and only consistency condition \eqref{eq:pia1} is satisfied, or, which is not surprising, if and only if the Perron -- Frobenius eigenvalue of the weighted admissibility matrix $\Uppi\Upsi$ is 1.
\end{rem}

\begin{rem}
Generally speaking, the base distribution of a multiplicative Markov measure need not be stationary with respect to the transition probabilities \eqref{eq:pab}. In the case of the free group and the associated topological Markov chain determined by the admissibility matrix $\Upsi$ \eqref{eq:SiF}, a description of the stationary measures has been essentially found in the course of the proof of \thmref{thm:pf}. Indeed, by \eqref{eq:coh}, a distribution $\th=\{\th_a\}_{a\in\Ac}$ on $\Ac$ is \textsf{$P$-stationary} (i.e., $\th P=\th$) if and only if $\th\Upnu^{-1}$ is a left Perron -- Frobenius eigenvector \eqref{eq:ph} of the operator $\Uppi\Upsi$, whence (up to a multiplicative constant)
$$
\th_a = \nu_a \ph_a = \nu_a \frac{1-\pi_{a^{-1}}}{1-\pi_a\pi_{a^{-1}}} \qquad\forall\,a\in\Ac  \;,
$$
or, in view of \eqref{eq:pia},
$$
\th_a = \nu_a (1 - \nu_{a^{-1}}) \qquad\forall\,a\in\Ac \;.
$$
In particular, the normalized stationary measure coincides with the base distribution $\nu=\{\nu_a\}$ if and only if the latter is uniform.
\end{rem}

\section{Harmonic measures of random walks}

\subsection{Boundary convergence}

A probability measure (\textsf{step distribution})~$\mu$ on a countable group $G$ determines the \textsf{random walk} $\RW(G,\mu)$ which is a Markov chain with the state space $G$ homogeneous both in time and space;  its transition probabilities are the translates of the measure $\mu$. In other words, the transitions
$$
g \mapstoto^{h\sim \mu} gh
$$
of this chain consist in the right multiplication by a $\mu$-distributed \textsf{increment}~$h$, so that the time $t$ position of a sample path issued at time 0 from the group identity~$e$ (the empty word) is the product
$$
g_t = h_1 h_2 \dots h_t
$$
of $t$ independent $\mu$-distributed increments $h_i$, and the distribution of $g_t$ is the $t$-fold \textsf{convolution power} $\mu^{*t}$ of the measure $\mu$. We denote by $\Pb$ the arising probability measure on the space of sample paths $\gb=(g_t)$ of the random walk $(\Fc,\mu)$ issued from the group identity~$e$ at time $0$, so that the time $t$ one-dimensional distribution of $\Pb$ is the $t$-fold convolution~$\mu^{*t}$. It is to this measure that we refer when talking about \emph{almost all} sample paths of the random walk.

We remind that a measure $\la$ on a $G$-space $X$ is called \textsf{$\mu$-stationary} if it is preserved by the convolution with $\mu$, i.e.,
$$
\la = \mu*\la = \sum_g \mu(g)\,g\la \;,
$$
where $g\la$ denotes the $g$-translate of the measure $\la$. Equivalently, it means that the measure $\la$ is stationary for the Markov chain on $X$ with the transitions $\displaystyle x \mapstoto^{h\sim\mu} hx$.

Returning to the free group $\Fc$, in this paper we are only interested in the situation when the random walk $(\Fc,\mu)$ is \textsf{nearest neighbour} (i.e., the support of the step distribution~$\mu$ is the set~$\Ac$ of free generators and their inverses). However, the very general argument below (we outline it for the sake of completeness) is applicable to \emph{all} step distributions $\mu$ that satisfy the following condition:
\begin{equation} \label{eq:mu}
\textit{the measure $\mu$ is not concentrated on a cyclic subgroup of $\Fc$}
\end{equation}
The cornerstone of this argument is a \textsf{contraction property} of the boundary action of the free group $\Fc$: for any sequence $g_n\in \Fc$ convergent to a point $\ga\in\p \Fc$ and any purely non-atomic probability measure $\th$ on $\p \Fc$ the sequence of translates $g_n\theta$ converges to the point measure $\de_\ga$ in the weak$^*$ topology.

\begin{thm} \label{thm:convergence}
If a step distribution $\mu$ on the free group $\Fc$ satisfies condition~\eqref{eq:mu}, then almost all sample paths of the random walk $(\Fc,\mu)$ converge to the boundary $\p \Fc$ of the group $\Fc$. Therefore, the sequence of the convolution powers $\mu^{*t}$ converges (weak$^*$) to the resulting hitting distribution (\textsf{harmonic measure}) $\nub$. The measure $\nub$ is purely non-atomic, and it is the unique $\mu$-stationary probability measure on $\p\Fc$.
\end{thm}

\begin{proof}[Proof (an outline)]
\phantom{z}

\begin{enumerate}[label = (\roman*)]
\item
The Krylov -- Bogolyubov theorem applied to the convolution operator $\la\mapsto\mu*\la$ implies that there exists a $\mu$-stationary probability measure $\la$ on $\p \Fc$.
\item
The set $A\subset\p\Fc$ of the maximal weight atoms of the measure $\la$ is finite and has the property that if $\ga\in A$, then also $h^{-1}\ga\in A$ for any $h\in\supp\mu$, whence $A$ is invariant with respect to the whole group $\gr\mu$ generated by $\supp\mu$.
\item
Condition \eqref{eq:mu} implies that the boundary action of the group $\gr\mu$ is \textsf{non-elementary} in the sense that there are no finite invariant subsets, and therefore the measure~$\la$ is purely non-atomic.
\item
For any function $\wh f\in C(\p \Fc)$ the corresponding \textsf{Poisson integral}
$$
f(g) = \left\lan \wh f, g\la \right\ran
$$
is a bounded $\mu$-harmonic function on $\Fc$. Thus, the sequence of its values
$$
f(g_n) = \left\lan \wh f, g_n\la \right\ran
$$
along the sample paths of the random walk $(\Fc,\mu)$ is a uniformly bounded martingale with respect to the natural coordinate filtration of the path space, and therefore it converges almost everywhere.
\item
Since the boundary $\p \Fc$ is separable, one concludes that for almost every sample path $\gb=(g_t)$ the sequence of translates $g_t\la$ converges in the weak$^*$ topology to a probability measure $\la_{\gb}$. Moreover, by passing to the limit as $t\to\infty$ in the identity
$$
\la = \mu^{*t}*\la = \int g_t\la \, d\Pb(\gb) \;,
$$
one obtains
\begin{equation} \label{eq:lambda}
\la = \int \la_{\gb}\,d\Pb(\gb) \;.
\end{equation}
\item
By (iii) and the aforementioned contracting property of the boundary action the convergence $g_t\la\to \la_{\gb}$ implies the convergence
\begin{equation} \label{eq:cg}
g_t\to g_\infty\in\p \Fc
\end{equation}
for a.e. sample path $\gb=(g_t)$, and
\begin{equation} \label{eq:delta}
\la_\gb=\de_{g_\infty} \;.
\end{equation}
\item
From \eqref{eq:cg} the distributions $\mu^{*t}$ of $g_t$ converge to the distribution $\nub$ of $g_\infty$. Furthermore, by \eqref{eq:lambda} and \eqref{eq:delta}
$$
\la = \int \de_{g_\infty}\, d\Pb(\gb) = \nub \;.
$$
\end{enumerate}
\end{proof}

\thmref{thm:convergence} is essentially due (in a somewhat implicit form though) to \textsc{Furstenberg} \cite[Corollary of Lemma~3.1]{Furstenberg71}, \cite[Section~14]{Furstenberg73}. Also see \textsc{Kaimanovich~-- Vershik} \cite[Section 6.8]{Kaimanovich-Vershik83}, \textsc{Cartwright~-- Soardi} \cite{Cartwright-Soardi89}, \textsc{Margulis} \cite[Chapter~6]{Margulis91}, \textsc{Woess} \cite{Woess93}, \textsc{Kaimanovich} \cite[Theorem~2.4]{Kaimanovich00a}, \textsc{Le\-drappier} \cite[Theorem 1.12]{Ledrappier01}. Actually, since the boundary $\p\Fc$ can also be interpreted as the \emph{space of ends} of the Cayley tree, for the nearest neighbour random walks that we consider in this paper the boundary convergence is an immediate consequence of the \emph{transience} of these random walks first established analytically by \textsc{Dynkin -- Malyutov} \cite{Dynkin-Malutov61}. The latter property is virtually obvious by today's standards and can be obtained in numerous ways (for instance, from the absence of invariant measures for the boundary action of the free group, see \textsc{Furstenberg} \cite[Section~9]{Furstenberg73}; of course, it also follows from \thmref{thm:convergence}). However, an advantage of the given proof of \thmref{thm:convergence} is not only in its generality, but also in providing as a byproduct uniqueness of the harmonic measure~$\nub$ as a $\mu$-stationary measure on the boundary $\p\Fc$. We use this fact in the proof of \thmref{thm:bij} below.

\subsection{Multiplicativity of the harmonic measure}

From now on we shall restrict ourselves to the situation when the support of the step distribution~$\mu$ on the free group~$\Fc$ is the set $\Ac$ of free generators and their inverses. Then $\mu$-stationarity of the harmonic measure $\nub$ implies its quasi-invariance with respect to the action of the group~$\Fc$, so that the Radon -- Nikodym \linebreak cocycle $\De_\nub$ of the harmonic measure is well-defined.

We denote by $\pi=\{\pi_a\}_{a\in\Ac}$ the collection of the \textsf{passage probabilities} of the random walk $(\Fc,\mu)$, i.e.,
\begin{equation} \label{eq:pass}
\pi_a = \Pb \bigl\{ \exists\, t>0: g_t=a \bigr\}
\end{equation}
is the probability that the random walk $(\Fc,\mu)$ issued from the group identity ever passes through (visits) $a\in \Ac$. The first observation is that the passage probabilities satisfy the equations
\begin{equation} \label{eq:pasys}
\pi_a = \mu_a + \left( \sum_{b\in\Ac\sm\{a\}} \mu_b \pi_{b^{-1}} \right) \pi_a \qquad \forall\,a\in\Ac
\end{equation}
(after the first $\mu$-distributed step either one is at $a$, or, in order to visit $a$ from $b\neq a$ one should first return to the group identity, and then go to $a$, cf. the proof of \thmref{thm:h} below). These equations imply that if one of the numbers $\pi_a$ is equal to~1, then all of them are 1, which would contradict transience of the random walk $(\Fc,\mu)$. Thus, $\pi_a<1$, and, on the other hand, $\pi_a\ge\mu_a>0$ for all $a\in\Ac$.

\begin{thm} \label{thm:h}
If $\mu=\{\mu_a\}_a\in\Ac$ is a probability measure on the free group~$\Fc$ with $\supp\mu=\Ac$, then the associated harmonic measure $\nub$ on $\p\Fc$ is multiplicative Markov (\dfnref{dfn:qp}), and its Radon -- Nikodym cocycle $\De_\nub$ is the product cocycle (\dfnref{dfn:product}) determined by the family of the passage probabilities \eqref{eq:pass}.
\end{thm}

\begin{proof}
We denote by $\nu=\{\nu_a\}_{a\in\Ac}$ the distribution of the first letter of the $\nub$-sampled infinite words, i.e., $\nu_a=\nu(\Sc_a)$ for $a\in\Ac$. Let us now consider the values of the harmonic measure on the shadows $\Sc_g$ of elements $g=a_1\dots a_n\in\Fc$ of length $n>1$. Since the increments of the random walk have length 1, if the limit point $g_\infty$ of a sample path $(g_t)$ belongs to $\Sc_g$, the path must have visited the point $a_1$. Therefore, the strong Markov property (e.g., see \textsc{Norris} \cite[Theorem 1.4.2]{Norris98}) applied to the first passage time
$$
T_{a_1} = \min \{t\ge 0: g_t=a_1 \}
$$
implies that
$$
\nu(\Sc_g) = \Pb \{ T_{a_1} < \infty \} \Pb_{a_1} \{ g_\infty \in \Sc_g \}
= \pi_{a_1} a_1\nu(\Sc_g) = \pi_{a_1} \nu (\Sc_{g'})
$$
for $g'= a_1^{-1} g = a_2 \dots a_n$, so that ultimately
\begin{equation} \label{eq:nu}
\nu(\Sc_g) = \pi_{a_1} \dots \pi_{a_{n-1}} \nu_{a_n} \;,
\end{equation}
whence the claim in view of \thmref{thm:final} and \thmref{thm:qi}.
\end{proof}

\begin{rem}
Alternatively, one can first identify the \emph{Martin boundary} of the random walk $(\Fc,\mu)$ with the geometric boundary $\p\Fc$, and then use the fact that the hitting distribution on the Martin boundary coincides with the representing measure of the constant harmonic function $\1$ as an integral of the minimal harmonic functions corresponding to boundary points. See \textsc{Dynkin~-- Malyutov} \cite[Theorem 3]{Dynkin-Malutov61} and \textsc{Dynkin -- Yushkevich} \cite[pp.~144-145, Problems 47, 48]{Dynkin-Yushkevich69} for a somewhat implicit description of the representing measure of an arbitrary positive $\mu$-harmonic function on $\Fc$ made explicit by \textsc{Levit -- Molchanov} \cite[formula (9)]{Levit-Molchanov71}, and the specialization of the latter to the case of the constant function $\1$ \cite[formula (10)]{Levit-Molchanov71}, as well as a modern exposition in \textsc{Woess} \cite[Chapter 9]{Woess09}.
\end{rem}

\begin{thm} \label{thm:bij}
Equations \eqref{eq:pasys} establish a bijective correspondence between probability measures $\mu=\{\mu_a\}_{a\in\Ac}$ with $\supp\mu=\Ac$ and collections of positive weights $\pi=\{\pi_a\}_{a\in\Ac}$ subject to conditions \eqref{eq:pia1} and \eqref{eq:sum}.
\end{thm}

\begin{proof}
As we have already noticed, the passage probabilities of the random walk $(\Fc,\mu)$ satisfy system \eqref{eq:pasys} and condition \eqref{eq:pia1}, whereas condition \eqref{eq:sum} follows from \thmref{thm:final} and \thmref{thm:qi}. Moreover, the passage probabilities provide the \emph{unique} solution of system \eqref{eq:pasys}. The reason is that in view of \dfnref{dfn:product} these equations can, after dividing by $\pi_a$, be rewritten as
$$
1 = \sum_{b\in\Ac}  \mu_a r_\pi(a,\ga) \qquad\forall\,\ga\in\p\Fc \;,
$$
where $r_\pi$ is the product cocycle determined by the family of weights $\pi$, which is precisely the $\mu$-stationarity condition for the solution of the Radon~-- Nikodym problem for the cocycle $r_\pi$, whereas such a solution exists and is unique by \thmref{thm:final} and \thmref{thm:qi}. Thus, uniqueness of the solution of system \eqref{eq:pasys} with respect to the variables $\pi_a$ for a fixed measure $\mu$ follows from the fact that the harmonic measure $\nub$ is the unique $\mu$-stationary probability measure on the boundary $\p\Fc$ (\thmref{thm:convergence}).

In the opposite direction, given a collection of weights $\{\pi_a\}_{a\in\Ac}$ that satisfy conditions \eqref{eq:pia1} and \eqref{eq:sum}, system \eqref{eq:pasys} is linear with respect to $\{\mu_a\}$ and can be easily solved explicitly. Indeed, one can express the weights $\mu_a$ from \eqref{eq:pasys} as
$$
\mu_a = \frac{(1-S)\pi_a}{1-\pi_a \pi_{a^{-1}}}
$$
by using the sum
$$
S=\sum_{a\in\Ac} \mu_a\pi_{a^{-1}} \;,
$$
for which then
$$
S = (1-S) Z \;,
$$
where
$$
Z = \sum_{a\in\Ac} \frac{\pi_a\pi_{a^{-1}}}{1-\pi_a\pi_{a^{-1}}} \;,
$$
so that finally
\begin{equation} \label{eq:mua}
\mu_a = \frac1{1+Z} \cdot \frac{\pi_a}{1-\pi_a\pi_{a^{-1}}} \qquad \forall\,a\in\Ac \;,
\end{equation}
and $\mu=\{\mu_a\}$ is a probability measure if and only if
$$
1 =  \sum_a \frac{\pi_a}{1-\pi_a\pi_{a^{-1}}} - Z = \sum_a \frac{\pi_a(1-\pi_{a^{-1}})}{1-\pi_a\pi_{a^{-1}}} \;,
$$
which is precisely condition \eqref{eq:sum}.
\end{proof}

\begin{rem}
The proof of \thmref{thm:bij} is based on the observation that, given a step distribution $\mu$ on $\Ac$, unique solvability of system \eqref{eq:pasys} for the weights $\{\pi_a\}$ is equivalent to existence and uniqueness of a $\mu$-stationary probability measure for the boundary action of the free group provided by \thmref{thm:convergence}. An analytic approach to solvability of the system \eqref{eq:pasys} goes back to \textsc{Dynkin -- Malyutov} \cite[Theorem~2]{Dynkin-Malutov61}, also see \textsc{Ledrappier} \cite[Lemma 2.2]{Ledrappier01} and \textsc{Mairesse} \cite[Theorem 4.5 and Lemma 4.7]{Mairesse05} (the \emph{traffic equations} system of Mairesse is precisely system \eqref{eq:pasys} expressed, by \eqref{eq:pia}, in terms of the base measure $\{\nu_a\}$ instead of the weights $\{\pi_a\}$).
\end{rem}

\begin{rem} \label{rem:cone}
The derivations of the expressions for the harmonic measure of shadow sets \eqref{eq:nu} and for the step distribution of the random walk in terms of the passage probabilities \eqref{eq:mua} in the proofs of \thmref{thm:h} and \thmref{thm:bij}, respectively, are applicable to any nearest neighbour random walk on an arbitrary tree (without any space homogeneity assumptions). As it was pointed out by an anonymous referee, the formula for recovering the transition probabilities from the first passage probabilities in this general setting already appears in \textsc{Cartier} \cite[formula (2.31)]{Cartier72}, and an explicit description of the harmonic measure was given in this situation by \linebreak \textsc{Nagnibeda -- Woess} \cite[Lemma 5.1]{Nagnibeda-Woess02}, also see a more detailed exposition in \textsc{Woess} \cite[Proposition~9.23 and formula (9.25)]{Woess09}.

Formally speaking, \emph{any} probability measure on the boundary of a tree is Markov for the chain whose sample paths are the geodesic rays issued from the tree root (the ``boundary process'' in the terminology of \cite{Nagnibeda-Woess02} and~\cite{Woess09}). However, the transition probabilities of the harmonic measure (either from \cite[formula (9.25)]{Woess09}, or from \cite[formula (5.2)]{Nagnibeda-Woess02}$\equiv$\cite[Corollary~9.57]{Woess09}; these two expressions are equivalent in view of the aforementioned formula of Cartier) are, in a sense, much more ``local''. More precisely, as stated in \cite[Theorem F$'$]{Nagnibeda-Woess02} or in \cite[p.~293]{Woess09}, they are entirely determined just by the corresponding \emph{cone types}. Then one can easily see that the harmonic measure is multiplicative in the sense that both its initial distribution and the transition probabilities are entirely determined by an appropriate base distribution on the set of cone types.
\end{rem}

\begin{rem} \label{rem:har}
The combination of \thmref{thm:bij} and \thmref{thm:final} immediately implies that for a product cocycle~$c$ solvability of the Radon -- Nikodym problem is equivalent to existence of a probability measure $\mu$ on the group~$\Fc$ (which can be chosen to be supported by the set $\Ac$, in which case $\mu$ is unique) such that the cocycle is \emph{$\mu$-harmonic}, i.e., $\sum_g \mu(g) c(g,\ga)=1$ for all boundary points $\ga\in\p\Fc$. One can call the cocycles with this property \emph{harmonic}.

Thus, for products cocycles solvability of the Radon -- Nikodym problem is equivalent to their harmonicity. One can easily deduce from a result of \textsc{Connell -- Muchnik} \cite[Theorem 1.1]{Connell-Muchnik07a} that this equivalence also holds for the cocycles of the boundary action whose potential is H\"older, and, even weaker, by \cite[Theorem 1.2]{Connell-Muchnik07a} this is also true for the cocycles which are continuously cohomologous to H\"older ones. It might be interesting to explore this relationship for other cocycles of the boundary action of the free group or for boundary actions of other groups (beginning, of course, with the word hyperbolic ones).

It is worth mentioning here that there is a superficially similar notion of harmonicity for cocycles with coefficients in unitary representations (e.g., see \textsc{Bekka} \cite{Bekka17} and the references therein), however conceptually it is quite different. In particular, our cocycles take values in the \emph{multiplicative} group of positive real numbers, whereas the cocycles from \cite{Bekka17} take values in the \emph{additive} group of the representation Hilbert space.
\end{rem}



\end{document}